\documentclass[11pt,reqno]{amsart}
\usepackage{amscd,graphicx,psfrag}
\usepackage{amsmath,amsfonts,epsfig}
\usepackage{url}

\newtheorem{defn}{Definition}[section]
\newtheorem{lemma}[defn]{Lemma}

\newtheorem{theorem}[defn]{Theorem}
\newtheorem{definition}[defn]{Definition}

\newtheorem{proposition}[defn]{Proposition}

\theoremstyle{definition}\newtheorem{example}[defn]{Example}

\newcommand{\zz}{{\mathbb Z}}

\newcommand{\nn}{{\mathbb N}}
\newcommand{\qq}{{\mathbb Q}}

\newcommand{\fn}{f^{(n)}}
\newcommand{\ff}{{\mathcal F}}
\newcommand{\Fhat}{\widehat F}
\newcommand{\Klift}{\widetilde{K}}

\newcommand{\D}{\mathcal D}

\newcommand{\Dlift}{\widetilde{\Delta}}

\newcommand{\ozsvath}{Ozsv\'{a}th}
\newcommand{\szabo}{Szab\'{o}}
\newcommand{\spin}{\ifmmode{\rm Spin}\else{${\rm spin}$\ }\fi}
\newcommand{\spinc}{\ifmmode{{\rm Spin}^c}\else{${\rm spin}^c$\ }\fi}
\newcommand{\spinct}{\mathfrak t}
\newcommand{\spincs}{\mathfrak s}

\newcommand{\bfa}{\mbox{\boldmath $\alpha$\unboldmath}}
\newcommand{\bfb}{\mbox{\boldmath $\beta$\unboldmath}}

\newcommand{\Ta}{{\mathbb T}_\alpha}
\newcommand{\Tb}{{\mathbb T}_\beta}

\newcommand{\bx}{\mathbf x}
\newcommand{\by}{\mathbf y}

\newcommand{\ds}{\displaystyle}

\DeclareMathOperator\Hom{Hom}

\DeclareMathOperator\cfhat{\widehat{CF}}
\DeclareMathOperator\hfhat{\widehat{HF}}
\DeclareMathOperator\cfkhat{\widehat{CFK}}
\DeclareMathOperator\hfkhat{\widehat{HFK}}
\DeclareMathOperator\cfk{CFK}
\DeclareMathOperator\hfk{HFK}

\DeclareMathOperator\hf{HF}

\newcommand{\qhs}{\mathbb{Q}\kern.3pt\mathrm{HS}^3}
\newcommand{\qhb}{\mathbb{Q}\kern.3pt\mathrm{HB}^4}

\newcommand{\intB}{\stackrel{\circ\,}{B^4}}

\newcommand{\lra}{\longrightarrow}
\newcommand{\lda}{\Big\downarrow}

\title[Knot concordance and branched covers]{Knot concordance and Heegaard Floer homology invariants in branched 
covers}
\author[J. Elisenda Grigsby]{J. Elisenda Grigsby}
\address{Department of Mathematics\newline\indent Columbia 
University\newline\indent 2990 Broadway MC4406\newline\indent NY, NY 
10027}
\email{\rm{egrigsby@math.columbia.edu}}
\author[Daniel Ruberman]{Daniel Ruberman}
\address{Department of Mathematics, MS 050 \newline\indent Brandeis 
University \newline\indent Waltham, MA 02454}
\email{\rm{ruberman@brandeis.edu}}
\author[Sa\v so  Strle]{Sa\v so  Strle}
\address{Faculty of Mathematics and Physics \newline\indent 
University of Ljubljana \newline\indent Jadranska 21 \newline\indent 
1000 Ljubljana, Slovenia }
\email{\rm{saso.strle@fmf.uni-lj.si}}
\subjclass[2000]{}
\date\today

\begin{document}
\begin{abstract}By studying the Heegaard Floer homology of the preimage of a knot $K \subset S^3$ inside its double branched cover, we develop simple obstructions to $K$ having finite order in the classical smooth concordance group.  As an application, we prove that all $2$-bridge knots of crossing number at most $12$ for which the smooth concordance order was previously unknown have infinite smooth concordance order.
\end{abstract}
\maketitle

\section{Introduction}

\bibliographystyle{halpha}

The introduction of a number of invariants has revitalized the study 
of classical knot concordance in recent years.  In this paper, we combine two of the most powerful of these new invariants, the correction term $d$ defined by \ozsvath-\szabo $\,$ in \cite{MR2222356} and the smooth concordance invariant $\tau$ defined by \ozsvath-\szabo $\,$ in 
\cite{MR2026543} and Rasmussen in \cite{GT0306378}, with older 
techniques due to Casson and Gordon~\cite{casson-gordon:orsay}, to 
give new obstructions to a knot being {\em smoothly} slice.   Our 
main results imply that many knots in the knot table 
(cf.~\cite{knotinfo,cha-livingston:unknown}) in fact have infinite 
order in the smooth concordance group.   Other than their common 
ancestry in the work of Casson-Gordon, our results do not overlap 
much with recent work on {\em topological} knot concordance stemming 
from the paper of Cochran-Orr-Teichner~\cite{MR1973052}; from now on 
we work exclusively in the smooth category.

Our approach relies on the study of the $d$ and $\tau$ invariants associated to the preimage of a knot $K \subset S^3$ in a cyclic branched cover (cf.~\cite{grigsby:hfkbranch,GT0610238}).    Let $K\subset S^3$, and $Y  \to S^3$ be the $n$-fold cyclic branched cover of $K$ with $n = q^r$ where $q$ is a prime.  It is well-known that $Y$ is a rational homology sphere ($\qhs$), and that if $K$ is slice, then the corresponding branched cover $W$ of the $4$-ball branched along the slice disk $\Delta$ is a rational homology ball ($\qhb$).  As observed in~\cite{MR1957829} and elaborated in~\cite{owens-strle:qhs4b,manolescu-owens:delta,jabuka-naik:dcover}, this implies that  many of the \ozsvath-\szabo\ $d$ invariants must vanish.

However, there is more information in the branched cover, as the 
branch set $\Klift$ is the boundary of the lift $\Dlift$ of $\Delta$ 
to $W$.   Now the knot $\Klift \subset Y$ has a collection 
$\{\tau_{\spincs}(Y,\Klift)\}$ of $\tau$ invariants, one for each 
$\spincs \in \spinc(Y)$.  In this situation, we prove the following vanishing 
theorem for those $\tau_{\spincs}$ associated to \spinc structures 
extending over $W$.  

\begin{theorem}\label{thm:slice-bound}
Let $K$ be a knot in $S^3$, and $Y$ the $q^r$-fold cover of $S^3$ 
branched along $K$. Denote by
$\Klift$ the preimage of $K$ in $Y$. If $K$ is slice, then there exists a 
subgroup $G<H^2(Y;\zz)$ with $|G|^2=|H^2(Y;\zz)|$
such that $d_{\spincs}(Y)=0$ and $\tau_{\spincs}(\Klift)=0$ for all 
$\spincs\in\spincs_0+G$, where
$\spincs_0$ is the spin structure on $Y$ described in Lemma~\ref{lem:spincover}.
\end{theorem}

As in many variations on the Casson-Gordon theme, it takes some work 
to extract computable obstructions from this theorem.  One important 
issue is that one does not know {\em a priori} which \spinc 
structures on $Y$ extend over $W$.    This problem gets worse when 
one studies the order of $K$ in the knot concordance group, because 
the number of \spinc structures that need to be examined can be very 
large.  (Compare the discussion in \S5 of \cite{jabuka-naik:dcover}.)  In addition, although there have been considerable advances in the computability of Heegaard Floer homology invariants in recent months (cf.~\cite{GT0607691},\cite{GT0607777}), much remains out of reach.

Even within the limited scope allowed by the current technology, however, one can deduce a considerable amount of information about knot concordance.  In particular, concentrating on the case $q^r = 2$, we define in Section \ref{sec:inforder}, for each prime $p$, two invariants, $\mathcal{T}_p(K) \in \mathbb{Z}$ and $\mathcal{D}_p(K) \in \mathbb{Q}$ which vanish on knots $K$ with finite smooth concordance order. These obstructions are particularly simple whenever $H^2(Y,\zz)$ is cyclic for $Y$ the double-branched cover of $K$ (as occurs, for example, for all $2$-bridge knots).  

\begin{theorem}\label{thm:obstdefn}
Let $K \subset S^3$ be a knot and $p \in \mathbb{Z}_+$ prime or $1$.  If there exists a positive $n \in \mathbb{Z}$ such that $\#_n K$ is smoothly slice, then $\mathcal{T}_p(K) = \mathcal{D}_p(K) = 0$.
\end{theorem}

Computation of these new invariants using the algorithm described in~\cite{GT0610238}  to compute $\tau$-invariants and the inductive formula in~\cite{MR1957829} for $d$-invariants allows us to determine the smooth concordance order of all $2$-bridge knots of $12$ or fewer crossings for which the smooth concordance order was previously unknown.

\begin{theorem}
All $2$-bridge knots\footnote{$K_{p,q}$ denotes the $2$-bridge knot whose double-branched cover is the lens space $-L(p,q)$.}  of $12$ or fewer crossings have smooth concordance order $1$,$2$ or $\infty$.
\end{theorem}

We are confident that as computational techniques improve, we will be able to gather similar results for a wider class of knots.  We remark that there are many methods (for example~\cite{jiang:cg,livingston-naik:4-torsion,livingston-naik:torsion,cochran-orr-teichner:structure}) for showing that a knot with finite order in the algebraic concordance group has infinite smooth or topological concordance order.   There is some overlap between results deduced by these methods and results from our paper.

The paper is organized as follows.  In Section \ref{sec:spin}, we discuss \spinc structures on branched covers of $B^4$.  In Section \ref{sec:vanish}, we give a proof of Theorem \ref{thm:slice-bound}, along with some relevant Heegaard Floer homology background.  In Section \ref{sec:inforder} we define $\mathcal{T}_p$ and $\mathcal{D}_p$ and prove that they provide obstructions to finite smooth concordance order.  We also discuss more general obstructions and comment briefly on the relationship of our results to recent work of Paolo Lisca \cite{lisca:ow}.  In Section \ref{sec:examples}, we apply our results to all $2$-bridge knots of $12$ or fewer crossings. We also show that the twist knots, with the exception of the figure-$8$ and Stevedore's knot, have infinite order in the knot
concordance group.  We conclude, in Section \ref{sec:comp2}, with details about the $\tau$ computations.

We thank Chuck Livingston for his help in the use of Knotinfo~\cite{knotinfo} and Jiajun Wang for some useful suggestions regarding the Floer homology computation algorithm.   The first author would also like to thank Peter \ozsvath $\,$ for some enlightening discussions during the course of this work.  The first author was partially supported by an NSF Postdoctoral Fellowship. The second author was partially supported by 
NSF Grant 0505605.   The third author was  supported in part by the Slovenian Research Agency program No. P1-0292-0101-04 and project No. J1-6128-0101-04.  Visits of the second and third authors were supported by a Slovenian-U.S.A. Research Project BI-US/06-07/003 and by the NSF.

\section{Lifting a spin structure}\label{sec:spin}
Theorem~\ref{thm:slice-bound} asserts the vanishing of the $\tau$ invariants on those \spinc structures on the $q^r$-fold covers of $S^3$ branched along a knot $K$ that extend over the corresponding branched cover of the $4$-ball.  The cohomology group  $H^2(Y;\zz)$ acts freely and transitively on the set of  \spinc structures on a manifold $Y$, and so a choice of one \spinc structure gives a bijection between that set and $H^2(Y;\zz)$.  We will show that there is a canonical \spinc structure, in fact a spin structure, on $Y$ that extends over the branched cover of the $4$-ball.  This is the spin structure $\spincs_0$ referred to in Theorem~\ref{thm:slice-bound}.
\begin{lemma}\label{lem:spincover}
Let $p:(W,\widetilde{F}) \to (B^4,F)$ be an $n$-fold cyclic branched cover with branch set a connected surface $F$.  Then there is a unique spin structure $\spincs_0$ on $W$ characterized as follows: if $n$ is odd, the restriction of $\spincs_0$ to $W-\nu(\widetilde{F})$ is the pullback $\widetilde\spincs$ of the spin structure on $B^4 -\nu(F)$ that extends over $B^4$, whereas if $n$ is even, the restriction of $\spincs_0$ to $W-\nu(\widetilde{F})$ is $\widetilde\spincs$ twisted by the element of $H^1(W-\nu(\widetilde{F});\zz_2)$ supported on the linking circle of $\widetilde{F}$. 
\end{lemma}

\begin{proof}
Let $\spincs$ be the spin structure on $B^4 -\nu(F)$ that extends over $B^4$ and let $\widetilde\spincs$ be its pullback to $W-\nu(\widetilde{F})$.
Then $\widetilde\spincs$ extends over $W$ iff its restriction to the circle bundle $S(\widetilde F)$ extends over the disk bundle $\nu(\widetilde F)$.
This happens precisely when the restriction to the fibre $S^1$ extends over the disk $B^2$.  Recall that a spin structure on an oriented manifold $M$ corresponds to a cohomology class in $H^1(F(M);\zz_2)$ where $F(M)$ is the oriented frame bundle.  Clearly $F(S^1)=S^1$, and one can check that the spin structure on $S^1$ that extends over $B^2$ corresponds to the nontrivial element of $H^1(S^1;\zz_2)$. Since the covering map is of order $n$, $\widetilde\spincs$ ``inherits'' this
property from $\spincs$ if $n$ is odd. If $n$ is even, $\widetilde\spincs$ does not extend over the disks. However, since the class of the fibre
is of infinite order in $H_1(W-\nu(\widetilde{F});\zz)$ we may twist the spin structure $\widetilde\spincs$ by the dual of this circle in 
$\Hom(H_1(W-\nu(\widetilde{F});\zz),\zz_2)$ to obtain a spin structure that extends over $W$.
\end{proof}

In the special case of a $2$-fold cover where the branch set is a disk, the cohomology $H^1(W;\zz_2)$ vanishes, so there is a unique spin structure on $W$.   The restriction of this spin structure to $Y$ is readily characterized as the unique invariant \spinc structure on $Y$, because the conjugation of \spinc structures acts on the odd-order group $H^2(Y;\zz)$ by multiplication by $-1$.   This observation is very convenient, because the methods of Section~\ref{sec:comp2} make it easy to identify this invariant \spinc structure.

\section{$\tau$ invariants of knots in rational homology spheres}\label{sec:vanish}
In this section we adapt the discussion from \cite{MR2026543} to nullhomologous knots in $\qhs$'s.  Our main aim is the proof of Theorem \ref{thm:slice-bound}.  We begin by collecting some standard Heegaard Floer homology facts.  See \cite{MR2113019}, \cite{MR2065507}, \cite{GT0512286} for more details.

Let $Y$ be a $\qhs$ and $K \subset Y$ an oriented 
null-homologous knot in $Y$.  We associate to the pair $(Y,K)$ a $2n$-pointed Heegaard diagram, i.e., a tuple $(\Sigma, \vec{\alpha}, \vec{\beta}, {\vec w}, {\vec z})$ associated to a handlebody decomposition of $Y$ coming from a generic self-indexing Morse function $f:Y \rightarrow \mathbb{R}$ with $|f^{-1}(0)| = |f^{-1}(3)| = n$.  Label the index $0$ critical points $a_1, \ldots a_n$ and the index $3$ critical points $b_1, \ldots b_n$.  Here, 
\begin{itemize}
\item
$\Sigma=f^{-1}(\frac{3}{2})$ is the Heegaard surface, 
\item
$\vec{\alpha} = (\alpha_1, \alpha_2, \ldots, \alpha_{g+(n-1)})$ are the coattaching circles of the $1$-handles, 
\item $\vec{\beta} = (\beta_1, \beta_2, \ldots, \beta_{g+n-1})$ are the attaching circles of the $2$-handles, and 
\item $\vec{w} = (w_1, \ldots, w_n)$ and $\vec{z} = (z_1, \ldots, z_n)$ are two $n$-tuples of points (all distinct) on $\Sigma - \vec{\alpha} - \vec{\beta}$, where $w_i$ specifies a unique flowline $\gamma_i$ from $b_i$ to $a_i$ (the one that intersects $\Sigma$ at $w_i$) and $z_i$ specifies a unique flowline $\eta_i$ from $b_{\sigma(i)}$ to $a_i$ ($\sigma$ some permutation of $\{1,\ldots, n\}$).
\end{itemize}

Then $K$ is uniquely determined by this data as the isotopy class of $${\bigcup}_{i=1}^n -\gamma_i \cup \eta_i.$$

In the case of a $2$-pointed Heegaard diagram, one produces the $\mathbb{Z} \oplus \mathbb{Z}$-filtered chain complex $\cfk^{\infty}(Y,K)$, a chain complex

\begin{itemize}
\item over\footnote{We define all Heegaard Floer chain complexes with $\mathbb{Z}_2$ coefficients in order to avoid orienting moduli spaces.} $\mathbb{Z}_2[U,U^{-1}]$,
\item whose generators are elements of the form $U^n{\bf x}$, where ${\bf x} \in \mathbb{T}_{\alpha} \cap \mathbb{T}_{\beta}$ is an intersection point between the $\alpha$ and $\beta$ tori in $Sym^g(\Sigma)$, and $n \in \mathbb{Z}$,
\item whose boundary map is given by 
$$\partial^{\infty}({\bf x}) = 
\sum_{{\bf y} \in \mathbb{T}_{\alpha} \cap \mathbb{T}_{\beta}}\sum_{\{\phi \in \pi_2({\bf x},{\bf y})|\mu(\phi)= 1\}} \#(\widehat{\mathcal{M}}(\phi))U^{n_w(\phi)}{\bf y}$$ 
where $\#(\widehat{\mathcal{M}}(\phi))$ is counted modulo $2$.\footnote{Our notation matches that of \cite{MR2065507}.}
\end{itemize}

Each intersection point ${\bf x}$ is assigned 

\begin{itemize}
\item an element ${\underline \spincs}({\bf x}) \in \spinc(Y_0(K))$ (where $Y_0(K)$ denotes $0$-surgery on $K$), 
\item $\spincs({\bf x}) \in \spinc(Y)$, 
\item an Alexander grading ${\bf A}({\bf x}) \in \mathbb{Z}$, and 
\item a Maslov grading ${\bf M}({\bf x}) \in \mathbb{Q}$.
\end{itemize}

Furthermore, the first three assignments are related by the natural splitting $$\spinc(Y_0) = \spinc(Y) \oplus \mathbb{Z}PD[\mu],$$ where $\mu$ is a choice of oriented meridian for $K$.  Thinking of a \spinc structure as a homology class of non-vanishing vector field (see \cite{MR1484699}, also Section 2.6 of \cite{MR2113019} and Section 2.3 of \cite{MR2065507}), the map $${\underline \spincs} \rightarrow \spincs$$ to the first summand is the unique extension to $Y$ of the restriction of ${\underline \spincs}$ to the knot complement, while the map $${\underline \spincs} \rightarrow {\bf A}$$ to the second is one-half of the evaluation of $c_1({\underline \spincs})$ on the surface obtained by capping-off a Seifert surface for $K$ with the core of the $0$-surgery.

The Maslov grading is the absolute $\mathbb{Q}$ homological grading defined in \cite{MR2222356}.

The $\mathbb{Z} \oplus \mathbb{Z}$-filtration on the $\cfk^{\infty}$ chain complex arises from an assignment of a filtration bigrading $A_w \times A_z$ to the chain complex generators (not to be confused with the ${\bf A} \times {\bf M}$ bigrading discussed in Section \ref{sec:comp2}).  This filtration bigrading is uniquely specified by the rules:

\begin{itemize}
\item $A_w({\bf x}) = 0 \,\, \forall \,\, {\bf x} \in \mathbb{T}_{\alpha} \cap \mathbb{T}_\beta$,
\item $A_z({\bf x}) = {\bf A}({\bf x}) \,\, \forall \,\, {\bf x} \in \mathbb{T}_\alpha \cap \mathbb{T}_\beta$,
\item $A_w(U^n{\bf x}) = A_w({\bf x}) - n$,
\item and $A_z(U^n{\bf x}) = A_z({\bf x}) - n$.
\end{itemize}

In addition, we remark that ${\bf M}(U^n{\bf x}) = {\bf M}({\bf x}) -2n$. The $E^2$ term of the associated spectral sequence is $\hfk^{\infty}(Y,K)$, and the $E^{\infty}$ term is $\hf^{\infty}(Y)$.

The $\cfhat(Y,K)$ chain complex is now the quotient complex corresponding to the $A_w =0$ slice of the $\cfk^{\infty}(Y,K)$ chain complex.  $\cfhat(Y,K)$ has a $\mathbb{Z}$-filtration coming from the $A_z = {\bf A}$ Alexander grading.  The associated graded chain complex of this filtered complex is $\cfkhat(Y,K)$.  Note that the $E^2$ term of the spectral sequence induced by the filtration on $\cfhat(Y,K)$ is $\hfkhat(Y,K)$, and the $E^{\infty}$ term is $\hfhat(Y)$.

Let $\ff(Y,K,\spincs,\ell)$ denote the subcomplex of $\cfhat(Y,K)$ generated by elements ${\bf x} \in \cfhat(Y;K)$ with $\spincs({\bf x})= \spincs$ and ${\bf A}({\bf x}) \leq \ell$. The inclusion $\ff(Y,K,\spincs,\ell)\to \cfhat(Y,\spincs)$, induces a morphism on homology denoted by $\imath_{K,\spincs}^{\ell}$. 

Recall that all of the aforementioned chain complexes (and, hence, all terms in the associated spectral sequences) split according to \spinc structures on $Y$: $$\cfk(Y,K) = \mathop{\oplus}\limits_{\{\spincs \in \spinc(Y)\}} \cfk(Y,K,\spincs).$$  Furthermore, $\hfhat(Y,\spincs)$ is isomorphic to $\zz_2$ if $Y$ is an L-space; for a general $\qhs$ $Y$, it contains a distinguished $\zz_2$ summand, which is in the image of $\hf^\infty(Y,\spincs)$.  We denote this copy of $\zz_2$ by $\hfhat_U(Y,\spincs)$.  

\begin{definition}\label{def:d}
The {\it correction term} for a torsion \spinc structure $\spincs$, denoted $d_{\spincs}(Y)$, is the absolute $\mathbb{Q}$ homological grading, {\bf M}, of $\hfhat_U(Y,\spincs)$.
\end{definition}

\begin{definition}\label{def:tau}
$\tau_{\spincs}(Y,K)=\tau_{\spincs}(K)$ is the minimal value of $\ell$ for which the image of $\imath_{K,\spincs}^{\ell}$ has nontrivial projection to $\hfhat_U(Y,\spincs)$.
\end{definition}

There is an interpretation of $\tau$ in terms of surgeries on $K$. Denote by 
$Q(K,\spincs,\ell)$ the
quotient complex $\cfhat(Y,\spincs)/\ff(Y,K,\spincs,\ell)$ and by 
$p_{K,\spincs}^\ell$
the map induced by the projection on homology. For any integer $n$ let
$$\Fhat_{n,\spincs,\ell}\colon 
\hfhat(Y,\spincs)\to\hfhat(Y_{-n}(K),\spincs_\ell)$$
denote the map associated to the two-handle cobordism, where the 
cobordism is endowed
with the \spinc structure $\spinct_{\ell}$, whose restriction to $Y$ 
is $\spincs$ and
which satisfies
$$\langle c_1(\spinct_\ell),[\widehat{S}]\rangle-n=2\ell\, ;$$
here $\widehat{S}$ denotes the surface obtained by capping-off a Seifert 
surface for $K$
with the core of the two-handle. These conditions uniquely specify 
$\spinct_\ell$ and hence
also the induced \spinc structure $\spincs_\ell$ on the surgery.

\begin{proposition}\label{prop:surgery-tau}
If $\ell<\tau_{\spincs}(K)$, then 
$\Fhat_{n,\spincs,\ell}(\hfhat_U(Y,\spincs))$ is nontrivial for all
sufficiently large $n$.
\end{proposition}

\proof
Fix a $2$-pointed Heegaard diagram for the pair $(Y,K)$ as well as a \spinc structure $\spincs \in \spinc(Y)$. Let $C_{\spincs}$ denote the $\cfk^{\infty}(Y,K,\spincs)$ chain complex.
Then 

\begin{itemize}
\item $C_{\spincs}\{A_w=0\}$ represents $\cfhat(Y,\spincs)$, and
\item $C_{\spincs}\{A_w=0,A_z \le\ell\}$ represents $\ff(Y,K,\spincs,\ell)$.
\end{itemize}

To relate this to the surgery recall Theorem 4.1 (and the discussion following it) in \cite{MR2065507}, which says 
that for all sufficiently large $n$, $\cfhat(Y_{-n}(K),\spincs_\ell)$ is identified with 
$C_{\spincs}\{\min(A_w,A_z-\ell)=0\}$.  Furthermore, under this identification the map $\Fhat_{n,\spincs,\ell}$ is 
induced by the
projection $f\colon C_{\spincs}\{A_w=0\}\to C_{\spincs}\{\min(A_w,A_z-\ell)=0\}$. This yields 
the following commutative
diagram
$$ \begin{array}{ccccccccc}
\vspace{2mm}
0 & \lra & C_{\spincs}\{A_w=0,A_z\le\ell\} & \lra & C_{\spincs}\{A_w=0\} & \lra & 
Q(K,\spincs,\ell) & \lra& 0 \cr
\vspace{2mm}
   &      &  \lda             &      & f \lda \    &      & \cong\, 
\lda  \     &     & \cr
0 & \lra & C_{\spincs}\{A_w\ge0,A_z=\ell\} & \lra & C_{\spincs}\{\min(A_w,A_z-\ell)=0\} & \lra & 
Q(K,\spincs,\ell) & \lra& 0 .
\end{array}
$$
Since for $\ell\le\tau_{\spincs}(K)$ the projection of the image of 
$\imath_{K,\spincs}^\ell$ into
$\hfhat_U(Y,\spincs)$ is trivial, $p_{K,\spincs}^\ell$ and therefore 
$\Fhat_{n,\spincs,\ell}$
are nontrivial on $\hfhat_U(Y,\spincs)$.
\endproof

The following properties of $\tau$ are important for our applications.

\begin{proposition}\label{prop:tau-properties}
(1) Let $(Y_i,K_i)$ be oriented knots and $\spincs_i\in\spinc(Y_i)$, $i=1,2$. Then
$\tau_{\spincs_1\#\spincs_2}(Y_1\# Y_2,K_1\# K_2)=\tau_{\spincs_1}(Y_1,K_1)+\tau_{\spincs_2}(Y_2,K_2).$\\
(2) If $(Y,K)$ is an oriented knot and $\spincs \in\spinc(Y)$, then
$\tau_{\spincs}(-Y,K)=-\tau_{\spincs}(Y,K)$.
\end{proposition}

\proof
(1) This follows from the K\"unneth Theorem for the knot filtration \cite[Theorem 7.1]{MR2065507} 
and from the fact that the tensor product of two vector space morphisms is nontrivial if and only if both of the
morphisms are nontrivial.\\
(2) If $(\Sigma,\bfa,\bfb,w,z)$ is a doubly pointed Heegaard diagram for $(Y,K)$, then
$(-\Sigma,\bfa,\bfb,w,z)$ is a diagram for $(-Y,K)$. If $\bx,\by\in\Ta\cap\Tb$ and
$\phi\in\pi_2(\bx,\by)$, then consider the homotopy class $\phi'$ of disks from $\by$ to $\bx$,
obtained from $\phi$ by precomposing with the complex conjugation. If $J_s$ is any one-parameter family
of complex structures, then $J_s$-holomorphic representatives for $\phi$ are in one-to-one
correspondence with $-J_s$-holomorphic representatives for $\phi'$. This correspondence induces
a duality map $\D\colon\cfhat_*(Y,\spincs)\to\cfhat^*(-Y,\spincs)$ which takes elements in $A_z$ filtration
level $\ell$ to $A_z$ filtration level $-\ell$; here $\cfhat_*$ denotes the chain complex and $\cfhat^*$ the cochain
complex obtained by applying the $\Hom(\cdot,\zz_2)$ functor. Applying the duality isomorphism to the
short exact sequence corresponding to the inclusion $\ff_*(Y,K,\spincs,\ell)\to\cfhat_*(Y,\spincs)$ yields
the following exact sequence 
$$0\to Q^*(-Y,K,\spincs,-\ell-1)\to\cfhat^*(-Y,\spincs)\to\ff^*(-Y,K,\spincs,-\ell-1)\to 0\, .$$
The image of the left map in homology has trivial projection into $\hfhat_U^*(-Y,\spincs)$ for all $\ell<\tau_{\spincs}(Y,K)$,
so $\ff_*(-Y,K,\spincs,-\ell-1)$ maps nontrivially to $\hfhat_{U,*}(-Y)$ for all $\ell$ satisfying $-\ell-1\ge-\tau_{\spincs}(Y,K)$. 
It follows that $\tau_{\spincs}(-Y,K)=-\tau_{\spincs}(Y,K)$.
\endproof

\begin{proposition}\label{prop:d-vanish}
Let $W$ be a $\qhb$ with boundary $Y$. Then
there exist $\spincs_0\in\spinc(Y)$ and a subgroup $G<H^2(Y;\zz)$ 
with $|G|^2=|H^2(Y;\zz)|$
such that for any $\spincs\in\spincs_0+G$ the following hold:\\
(1) $d_{\spincs}(Y)=0$ and\\
(2) $\spincs$ extends to $\spinct\in\spinc(W)$ and the map 
$\Fhat_{W-\intB,\spinct}\colon
\hfhat(S^3)\to\hfhat_U(Y,\spincs)$ is an isomorphism.\\
Moreover, if $W$ is spin, then $\spincs_0$ can be chosen to be a spin 
structure.
\end{proposition}

\proof
The first statement follows from \cite[Proposition 
4.1]{owens-strle:qhs4b}. For the second
note that $F^\infty_{W-\intB,\spinct}\colon 
\hf^\infty(S^3)\to\hf^\infty(Y,\spincs)$ is an
isomorphism according to \cite[Theorem 9.6]{MR2222356}. Since $W$ is 
a $\qhb$,
this map preserves absolute grading, hence the map in degree 
zero is an isomorphism.
\endproof

\begin{theorem}\label{thm:tau-bound}
Let $W$ be a $\qhb$ with boundary $Y$ and 
$S\subset W$ a surface
whose boundary is a knot $K\subset Y$. If $\spincs\in\spinc(Y)$ 
extends over $W$ and
$d_{\spincs}(Y)=0$, then $g(S)\ge|\tau_{\spincs}(K)|.$
\end{theorem}

\proof
If $g(S)=0$ we may replace $K$ with the connected sum of $K$ and 
the trefoil of appropriate
handedness so as to increase both sides of the claimed inequality by 
1. Thus we may assume
$g(S)>0$. We may further assume that $\tau_{\spincs}(Y,K)\ge0$ by
changing the orientation of $Y$ if necessary (see Proposition \ref{prop:tau-properties}).

Choose $\ell<\tau_{\spincs}(K)$ and let $X$ be obtained from $W$ by 
adding a two-handle along $K$
with framing $-n$, where $n$ is large enough so that Proposition 
\ref{prop:surgery-tau} applies.
Endow $X$ with a \spinc structure $\spinct$, whose restriction to $Y$ 
is $\spincs$ and whose
restriction to the two-handle cobordism agrees with $\spinct_\ell$ 
(defined before Proposition
\ref{prop:surgery-tau}). Then
$$\langle c_1(\spinct),[\widehat{S}]\rangle-n=2\ell\, ,$$
where $\widehat{S}$ denotes the surface obtained by capping-off 
$S$ with the core of the
two-handle. The map $\Fhat_{X-\intB,\spinct}\colon 
\hfhat(S^3)\to\hfhat(Y_{-n}(K),\spincs_\ell)$
is nontrivial, since it is the composition of the map induced by $W$ with
image $\hfhat_U(Y,\spincs)$ and the map induced by the two-handle cobordism that
is nontrivial on $\hfhat_U(Y,\spincs)$ by Proposition \ref{prop:surgery-tau}.

Now we split the cobordism $X$ differently. 
Let $N$ be a tubular neighborhood of $\widehat{S}$. Then by 
\cite[Lemma 3.5]{MR2026543} the map
induced by the cobordism $N-\intB$ and \spinc structure $\spinct$ is 
trivial unless
$$\langle c_1(\spinct),[\widehat{S}]\rangle-n\le 2g(S)-2\, ,$$
from which we conclude $|\tau_{\spincs}(K)|\le g(S)$. 
\endproof

We note that Theorem \ref{thm:slice-bound} is an immediate consequence
of Proposition \ref{prop:d-vanish} and Theorem \ref{thm:tau-bound}.

\section{Obstructions to finite concordance 
order}\label{sec:inforder}

Using Theorem~\ref{thm:slice-bound} to test 
whether a given knot is slice might in principle require a good deal of 
calculation. This is because one does not know in advance the 
subgroup $G$ with $|G|^2 = |H^2(Y;\zz)|$ (referenced in the theorem) on which all $\tau$ and $d$ invariants must vanish; hence, to rule out the existence of such a subgroup, one must find all subgroups $G$ of the appropriate order and verify that, indeed, either $\tau_{\spincs}(\widetilde{K}), d_{\spincs}(Y) \neq 0$ for some $\spincs \in (\spincs_0 + G)$ in each case.  

This can be computationally formidable; as an example 
consider the $2$-bridge knot $K_{45,17}$.  The $d$-invariant 
calculations in~\cite{jabuka-naik:dcover} left open the possibility 
that $K = \#_4 K_{45,17}$ might be 
slice.  Similarly, the $\tau$ obstruction for $K$ necessitated the computation of $\tau_{\spincs}$ for $\spincs \in \spincs_0 + G$ on all order $45^2$ subgroups $G$ of $(\zz/45)^4$, which 
we did using the program pari-gp~\cite{PARI2}.  The computation took 
quite a while (there are $9,745,346$ such subgroups), and indeed there exist 
nonzero $\tau$-invariants in each one so that $K$ is not slice. 
The complexity of such calculations clearly gets out of hand rapidly.

Inspired by Jabuka-Naik's idea (see Obstruction 5.1 in \cite{jabuka-naik:dcover}) of looking at smaller-order, more computationally-accessible, subgroups of $H^2(Y;\mathbb{Z})$, we developed two simple obstructions to a knot having finite smooth concordance order.  Although at first glance these obstructions seem much weaker than the original obstructions provided by Theorem \ref{thm:slice-bound}, they have had remarkable success on the class of knots upon which we were able to perform calculations.  In fact, our obstructions were able to show that all of the $2$-bridge knots considered in~ \cite{jabuka-naik:dcover} have infinite concordance order.  We also determined the (previously unknown) smooth concordance order of all $2$-bridge knots with crossing number at most 12.  The obstructions and examples follow in the next two sections.  The final section is devoted to a brief explanation of how the $\tau$ calculations were performed.

We will find the following notation useful in what follows. If $f\colon A \to \qq$ is a function on a finite abelian group and $H<A$ is any subgroup
we let $S_H(f)=\sum_{h\in H}f(h)$.

\begin{definition}\label{def:TD}
Let $K \subset S^3$ be a knot, $Y$ the double-branched cover of $K$, $\widetilde{K}$ the preimage of $K$ in $Y$, and $p \in \mathbb{Z}_+$ either a prime or $1$.  Fix an affine identification of $\spinc(Y)$ with $A = H^2(Y;\zz)$ such that the distinguished spin structure $\spincs_0$ mentioned in Lemma \ref{lem:spincover} corresponds to 0.  

Let $\mathcal{G}_p$ denote the set of all order-$p$ subgroups of $A$. Define
$$\mathcal{T}_p(K) := \left\{\begin{array}{cl}
                                   \ds{\min\left\{\left|\sum_{H \in \mathcal{G}_p} n_H S_H(\tau(Y, \widetilde{K}))\right|\, ;\, n_H \in\zz_{\ge 0}\right\}} & \mbox{if p divides det(K)}\\
                                   0 & \mbox{otherwise}
                                 \end{array}\right.$$
and

$$\mathcal{D}_p(K) := \left\{\begin{array}{cl}
                                  \ds{\min\left\{\left|\sum_{H \in \mathcal{G}_p} n_H S_H(d(Y))\right|\, ;\, n_H \in\zz_{\ge 0}\right\}} & 
                                   \mbox{if p divides det(K)}\\
                                   0 & \mbox{otherwise}
                                 \end{array}\right.$$

\end{definition}

It is worthwhile to remark at this point that the definitions of $\mathcal{T}_p$ and $\mathcal{D}_p$ are considerably simpler when $A = H^2(Y;\zz)$ is cyclic (as is the case for all $2$-bridge knots).  In this case, there is a unique subgroup of $A$ of order $p$ for each $p$ dividing $det(K)$, and hence $\mathcal{T}_p$ (resp. $\mathcal{D}_p$) is just the absolute value of 
$$\sum_{\{\spincs \in \spinc(Y)| \spincs \,\, \mbox{{\tiny has order p}}\}} \tau_{\spincs}(\widetilde{K}) \,\qquad \mbox{(resp.} \,\, d_{\spincs}(Y)\mbox{)}.$$
More generally we note that $\mathcal{T}_p$ or $\mathcal{D}_p$ is nonzero if the sums $S_H$ of the corresponding invariants are nonzero and of the same sign on all the subgroups $H$ of order $p$.

Note that the invariant $\mathcal{D}_1$ is $\frac{1}{2}\delta$, where $\delta$ is Manolescu-Owens' concordance invariant \cite{manolescu-owens:delta}.

The rest of this section will be devoted to proving Theorem \ref{thm:obstdefn}, the statement of which we repeat here for convenience:

{\flushleft {\bf Theorem \ref{thm:obstdefn}.}  {\it  Let $K \subset S^3$ be a knot and $p \in \mathbb{Z}_+$ prime or $1$.  If there exists a positive $n \in \mathbb{Z}$ such that $\#_n K$ is smoothly slice, then $\mathcal{T}_p(K) = \mathcal{D}_p(K) = 0$.}} 
\begin{proof}
The proof is essentially Theorem \ref{thm:slice-bound} in the case $q^r=2$ combined with the elementary observation that a finite abelian group of order $m$ contains a subgroup of order $p$ ($p$ prime) whenever $p|m$.  

Assume that $K$ has finite smooth concordance order, $n$, and fix a prime $p$.

Let $A$ denote $H^2(Y;\zz)$, where $Y$ is the double-branched cover of $K$.  By Theorem \ref{thm:slice-bound} in the case $q^r = 2$, there exists a subgroup $G < A^n$ with $|G| = |A|^{{n/2}}$ on which $\tau$ and $d$ vanish identically (here, again, we have fixed an affine identification of $\spinc(Y)$ with $H^2(Y;\zz)$).  As usual, $\widetilde{K}$ denotes the preimage of $K$ in $Y$.

We represent an element $g \in G$ by $g = (g_1, \ldots, g_n)$ where $g_i \in A$.  Note that $$\tau_g(\#_n\widetilde{K}) = \sum_{i=1}^n \tau_{g_i}(\widetilde{K})$$ and $$d_g(\#_n Y) = \sum_{i=1}^n d_{g_i}(Y)$$ (see Proposition \ref{prop:tau-properties} and Section 2.3 in \cite{jabuka-naik:dcover}).

Consider a finite abelian group, $A$, and suppose that we have a function $f:A \to \qq$.  We have in mind $A = H^2(Y;\zz)$ and $f$ either the $\tau$ or $d$ invariant subject to a fixed affine identification of $\spinc(Y)$ with $A$.  Given such an identification, there is a straightforward extension that identifies $\spinc(\#_n Y)$ with $A^n$.  For $n \in \nn$, denote by $\fn: A^n \to 
\qq$ the function $\fn(g_1,\ldots,g_n) = f(g_1)+ \cdots f(g_n)$. 

We assume without loss of generality that $p$ divides $det(K)$ and hence $|G|$.  Therefore, $G$ contains an element of order $p$, say $g = (g_1, \ldots g_n).$  Note that $\forall \,\, i \in \{1, \ldots n\}$ $g_i$ has order dividing $p$ and at least one of the $g_i$ has order $p$.  Let $<g> = \{0g, g, \ldots, (p-1)g\}$ denote the cyclic subgroup in $G$ generated by $g$ and $G_i$ denote the cyclic subgroup in $A$ generated by $g_i$.  Theorem \ref{thm:slice-bound} tells us that $\fn$ vanishes identically on $G$, hence on $<g>$ for $f=\tau,d$.  In particular, $$\fn(mg)=0 \,\, \mbox{for} \,\,m=0,\ldots, p-1.$$
Note that $\fn(0g) = nf(0) = 0$ implies that $\mathcal{T}_1 = \mathcal{D}_1 = 0$, proving the proposition for $p=1$.

We now have:
$$
\begin{array}{rcll}
\vspace{10pt}
  \fn(mg) &=& 0 \,\,\forall\,\, m = 0, \ldots, p-1 \,\, & \Longrightarrow\\
\vspace{10pt}
  \sum_{m=0}^{p-1} \fn(mg) &=& 0 \,\, & \Longrightarrow\\
\vspace{10pt}
  \sum_{m=0}^{p-1}\sum_{i=1}^{n}f(mg_i) &=& 0 \,\, & \Longrightarrow\\
\vspace{10pt}
  \sum_{i=1}^n\sum_{m=0}^{p-1}f(mg_i) &=& \sum_{i=1}^n S_{G_i}(f)=0 & \\
\end{array}
$$
\end{proof}

As in Jiang's proof~\cite{jiang:cg} that algebraically slice knots form an infinitely generated subgroup of the concordance group, this test can be applied, one prime at a time, to show that knots are linearly independent.  See Proposition~\ref{prop:twist-knots} for an application of this principle.

\subsection{Further obstructions to finite concordance order}\label{sec:furthobst}
Even when $\mathcal{T}_p$ and $\mathcal{D}_p$ vanish, it is sometimes possible, through more careful analysis, to find an obstruction to finite concordance order.  The following two propositions describe tests we developed to deal with the knots $K_{77,18}$, $K_{81,14}$, $K_{125,33}$, and $K_{209,81}$, for which $\mathcal{T}_p$ and $\mathcal{D}_p$ failed to provide an obstruction.  

\begin{proposition}\label{prop:minmax}
Let $K \subset S^3$ be a knot of finite concordance order and let $p$ be a prime. Denote by $Y$ the 
double-branched cover of $K$ and by $\Klift$ the lift of $K$ to $Y$. Suppose that the $p$-subgroup 
$A_p$ of $A=H^2(Y;\zz)$ is isomorphic to $\zz_p$ and fix some affine identification of $\spinc(Y)$ 
with $A$ that sends the spin structure on $Y$ to $0$.  Then
$$\min\{f(\spincs)\, ; \, \spincs\in A_p \} = -\max\{f(\spincs)\, ; \, \spincs\in A_p \},$$
where $f$ denotes either $\tau(\Klift)$ or $d(Y)$.
 
Moreover, let $\Delta(\pm M)=\{d\in\zz_p\, ;\, d\ne 0,\ f(a)=f(a+d)=\pm M \text{\ for some $a\in A_p$}\}$, 
where $M$ denotes the maximum of $f$ on $A_p$. Then 
$$\bigcap_{d\in\Delta(M)} d^*\Delta(-M)\ne \emptyset\, ,$$
where $d^*$ denotes the (multiplicative) inverse of $d$ modulo $p$.
\end{proposition}
\begin{proof}
Suppose $\#_{2n}K$ is smoothly slice. 
If $G<A^{2n}$ has order $|A|^n$, then $G$ contains a subgroup $G_p$ isomorphic to $(\zz_p)^n$.
Note that $G_p$ is a subgroup of $H=(\zz_p)^{2n}<A^{2n}$.
Let $\{g_i\, ;\, i=1,\ldots, n\}$ be a set of generators for $G_p$. By elementary operations
and rearrangement of summands in $H$ we may assume that the generators are of the form $g_i=(e_i,h_i)$
where $e_i\in(\zz_p)^n$ has the only nonzero entry in the $i^{\text{th}}$ component equal to $1$
and $h_i=(h_{ij})\in(\zz_p)^n$. 

Define $M$ (resp.\ $m$)
to be the maximum (resp.\ minimum) of $f$ on $A_p$. Assume contrary to the statement of the proposition that 
$M+m\ne 0$. Then by replacing $f$ with $-f$ if necessary, we may assume that $M+m>0$.
Let $k\in\zz_p$ be such that $f(k)=M$. Then
$$f^{(2n)}\Big(\,\textstyle{\sum\limits_{i=1}^nkg_i}\Big)=nM+\ds{\sum_{j=1}^n}f\Big(\,\textstyle{\sum\limits_{i=1}^nkh_{ij}}\Big)>0,$$
which contradicts Theorem \ref{thm:slice-bound}.

Fix $\ell\in\{1,\ldots,n\}$, $k\in\zz_p$ with $f(k)=M$, and choose some $d\in\Delta(M)$. Since 

\begin{eqnarray*}
f^{(2n)}\Big(\,\textstyle{\sum\limits_{i=1}^nkg_i}\Big)&=&f^{(2n)}\Big(dg_\ell+\textstyle{\sum\limits_{i=1}^nkg_i}\Big)\\
&=& nM + \sum_{j=1}^n f(a_j)\\
&=&0
\end{eqnarray*}
implies that $f(a_j) = -M$ for all $j$ (where $a_j = k(\sum_{i=1}^n h_{ij})$), it follows that $dh_{\ell j}\in\Delta(-M)\cup\{0\}$ for $j=1,\ldots,n$. Since $g_{\ell}$ is nonzero, at least one $h_{\ell j}$ has to be
nonzero and hence $h_{\ell j}\in d^*\Delta(-M)$.  In fact, since the above did not depend on which $d \in \Delta(M)$ we chose, we conclude that $$h_{\ell j} \in \bigcap_{d \in \Delta(M)} d^*(\Delta(-M)),$$ as desired.
\end{proof}

In another direction, we can extend the definitions of $\mathcal{T}_p$ and $\mathcal{D}_p$ to include subgroups of prime-power order.  More specifically, assuming notation from Definition \ref{def:TD}, suppose that $p^k$ divides the exponent of $A$ for some $k>1$.
Then letting $\mathcal{G}_{p^k}$ be the set of all cyclic subgroups of $A$ of order $p^k$
one may define $\mathcal{T}_{p^k}$ and $\mathcal{D}_{p^k}$ as in the case $k=1$. The conclusion of
Theorem \ref{thm:obstdefn} can then be strengthened to $\mathcal{T}_{p^k}(K) = \mathcal{D}_{p^k}(K) = 0$
for all $k$ for which $p^k$ divides the exponent of any subgroup of $H^2(Y;\zz)^n$ of order
$|H^2(Y;\zz)|^{n/2}$ (as usual, $Y$ is the double-branched cover of $K$). The largest power of $k$
as above is difficult to determine in general; we state a proposition for the particular class of groups of interest to us here.

\begin{proposition}\label{prop:primepwr}
Let $K \subset S^3$ be a knot and $p$ a prime. Suppose that the $p$-subgroup of $H^2(Y;\zz)$ is cyclic 
and let $p^m$ be the largest power of $p$ dividing $|H^2(Y;\zz)|$. If some positive multiple of $K$ is 
smoothly slice, then $\mathcal{T}_{p^k}(K) = \mathcal{D}_{p^k}(K) = 0$ for all $k\le \lfloor \frac{m+1}2\rfloor$.
\end{proposition}
\begin{proof}
Suppose $A$ is a finite abelian group with cyclic $p$-subgroup and $|A|=p^ms$, $(p,s)=1$. If $G<A^n$ has
order $|A|^{n/2}$, then $p^{mn/2}$ divides the order of $G$. It follows that either $G$ has an element of
order $p^{m/2}$ (if $m$ is even) or of order $p^{(m+1)/2}$ (if $m$ is odd). Clearly then one of the
components of this element has that order.

From here the proof proceeds inductively on $k$ where each step is essentially the same as the case
$k=1$; the only difference being that for $k>1$ at the end of the proof we get a sum of sums $S_H$, 
where the order of $H$ divides $p^k$.
\end{proof}

Although one may be tempted to extend the definitions of $\mathcal{T}_p$ and $\mathcal{D}_p$ to more general $p$, there are problems with doing so.  As an example, consider the $2$-bridge knot $K_{45,17}$. Its $d$ invariants vanish on the order $3$ subgroup and have sum zero on the order $5$ subgroup, but the sum of the $d$ invariants over any larger subgroup is nonzero. One may wonder if that could be used as an obstruction to finite concordance order, since any subgroup of order $45^n$ in $(\zz_{45})^{2n}$ contains either elements of order 9 or elements of order 15. Using an analogue $\mathcal{D}_{15}$
(and $\mathcal{D}_9$) of the invariant $\mathcal{D}$ from Definition \ref{def:TD} one would like to conclude that indeed the knot has infinite concordance order. Unfortunately, there are elements 
of order 15 in $(\zz_{45})^{2n}$ that have components of orders 3 and 5 only. In fact, with this in mind, it is easy to construct subgroups of order $45^2$ in $(\zz_{45})^4$ 
on which the obstruction coming from $d$ invariants vanishes. This is consistent with results of 
\cite{jabuka-naik:dcover}.

\subsection{Lisca's work}
It is worth commenting upon the relationship between our results and recent remarkable work of Paolo Lisca.  He shows in ~\cite{lisca:ow} that if $K_{p,q}$ is slice, then it is on a list of knots previously known to be ribbon.  He uses the fact that if $K_{p,q}$ is slice, then its double cover $-L(p,q)$ bounds a rational homology ball.  By Donaldson's theorem, the intersection form, $P$, of a canonical plumbed $4$-manifold bounded by $-L(p,q)$, embeds in a diagonal form, $D_n$.   A complicated algebraic argument then produces the list of possible $(p,q)$.  Lisca's theorem does not {\em a priori} decide the concordance order of any $2$-bridge knot, and so our results do not follow from his.  It is possible that a generalization of his method would say something about concordance order, but the following example suggests that this may be difficult to carry out; it shows that the assumption that $P \oplus P$ embeds in $D_n\oplus D_n$ does not imply that $P$ embeds in $D_n$.

\begin{example} Let $P$ be the intersection form of the plumbing of a degree $2$ and a degree $3$ bundle over $S^2$. Then $P$ does not embed in $D_2$, for you need at least $3$ vectors to realize intersection number $3$.  However, $P \oplus P$ can be embedded in $D_2 \oplus D_2=D_4$ as follows: let $e_1,\ e_2,\ e_3,$ and $e_4$ denote the standard ON basis and let $a_i$ and $b_i$ be generators in $P_i=P$ for $i=1,\ 2$ with $a_i^2=2$ and $b_i^2=3$. Then the map
\begin{align*}a_1 \to & e_1+e_2\\b_1 \to  & e_1+e_3+e_4\\a_2 \to &  e_3-e_4\\b_2 \to  & -e_1+e_2+e_3\\\end{align*}
gives an embedding of $P \oplus P$ in $D_4$.
\end{example}

\section{Examples} \label{sec:examples}
There are many knots in the knot tables whose order in the smooth 
concordance group is unknown~\cite{knotinfo,cha-livingston:unknown}.  Combining Theorem \ref{thm:slice-bound} with 
\begin{itemize}
\item calculations of $\tau_{\spincs}(\widetilde{K})$ for $\spincs \in \spinc(Y)$ ($Y$ the double-branched cover of a $2$-bridge knot $K$) (described in 
Section~\ref{sec:comp2}) and
\item previously-known calculations of $d_{\spincs}(Y)$ ($Y$ again the double-branched cover of a $2$-bridge knot $K$) using the inductive formula in Section 4.1 of \cite{MR1957829}, 
\end{itemize}
\begin{table}[ht]\caption{}
\label{table:tableofknots}\begin{tabular}{||c|c|c|l||} \hline Knot $K$
&$2$-bridge notation  &  Order of $K$ & Test \cr \hline \hline 
$8_{13}$ & $29/11$ & $\infty$  & $\mathcal{T}_{29} \neq 0$  and $\mathcal{D}_{29} \neq 0$ \cr 
\hline$9_{14}$ & $37/14$ &  $ \infty$ & $\mathcal{T}_{37} \neq 0$ and $\mathcal{D}_{37} \neq 0$ \cr 
\hline$9_{19}$ & $41/16$ & $ \infty$ & $\mathcal{T}_{41} \neq 0$ and $\mathcal{D}_{41} \neq 0$ \cr 
\hline$10_{10}$ & $45/17$ & $ \infty$ & $\mathcal{T}_3 \neq 0$ and $\mathcal{T}_5 \neq 0$ \cr 
  &   &  &but $\mathcal{D}_3 = 0$ and $\mathcal{D}_5 = 0$ \cr 
\hline$10_{13}$ & $53/22$ & $ \infty$ & $\mathcal{T}_{53} = 0$, but $\mathcal{D}_{53} \neq 0$ \cr
\hline$10_{26}$ & $61/17$ & $\infty$ &  $\mathcal{T}_{61} \neq 0$  and $\mathcal{D}_{61} \neq 0$\cr 
\hline$10_{28}$ & $53/19$ & $\infty$ & $\mathcal{T}_{53} \neq 0$   and $\mathcal{D}_{53} \neq 0$\cr 
\hline$10_{34}$ & $37/13$ & $\infty$ & $\mathcal{T}_{37} \neq 0$  and $\mathcal{D}_{37} \neq 0$\cr 
\hline$11_{91}$ & $129/50$ & $\infty$ & $\mathcal{T}_3 \neq 0$  and $\mathcal{D}_{3} \neq 0$ \cr 
\hline$11_{93}$ & $93/41$ & $ \infty$ & $\mathcal{T}_3 \neq 0$  and $\mathcal{D}_{3} \neq 0$\cr 
\hline$11_{98}$ & $77/18$ & $\infty$ & Proposition~\ref{prop:minmax} for $p = 7,\; 11$\cr \hline 
$11_{119}$ & $77/34$ & $ \infty$ & $\mathcal{T}_{11}  \neq 0$ and $\mathcal{D}_{11}\neq 0$\cr 
\hline\end{tabular}\end{table}
we have been able to determine all of the previously-unknown smooth concordance orders for $2$-bridge knots
of at most $12$ crossings.  Some sample results are summarized in Table \ref{table:tableofknots}, where we also indicate which invariants provide an obstruction to finite smooth concordance order.   

The most interesting example from our point of view was $K = 10_{10} = K_{45,17}$.  As remarked above (Subsection~\ref{sec:furthobst}) the $d$-invariants fail to obstruct the possibility that $K$ has order $4$.    However, $\tau(\Klift)$ shows $K$ has infinite concordance order, since both $\mathcal{T}_3 \neq 0$ and $\mathcal{T}_5 \neq 0$.    Further along in the Knotinfo tables, we found the $12$-crossing $2$-bridge knots $K_{81,14}$ and $K_{125,33}$, for which $\mathcal{D}_3$ and $\mathcal{D}_5$ vanish (respectively).  However, $\mathcal{D}_9$ and $\mathcal{D}_{25}$ are non-zero, and using Proposition \ref{prop:primepwr} we showed that these knots have infinite order. Similarly for the $12$-crossing $2$-bridge knot $K_{209,81}$ both $\mathcal{T}$ and $\mathcal{D}$ invariants associated to $11$ and $19$ are zero, however the knot has infinite order by Proposition \ref{prop:minmax} applied to either $d$ or $\tau$ invariants.   We remark that some of the knots that we treated (for instance $K_{77,18}$ and $K_{209,81}$ can be shown to have infinite topological concordance order by using the main result of Livingston-Naik~\cite{livingston-naik:4-torsion}.

\subsection{Twist Knots}\label{sec:twist}
We mention the subclass of twist knots ($K_{p,2}$), since these were the $2$-bridge knots originally addressed by Casson and Gordon.
The following result is a generalization of Jiang's theorem \cite{jiang:cg} that the set of algebraically slice twist knots $K_{p^2,2}$, $p\ge 5$ a prime, is
linearly independent in the concordance group.

\begin{proposition}\label{prop:twist-knots}
All twist knots $K_{p,2}$, $p\ge 3$, have infinite order in the knot concordance group except for $p=9$ (Stevedore's knot, which is slice)
and $p=5$ (figure-$8$ knot, which has order $2$). Moreover, twist knots in any family $K_{p_i,2}$ with $p_i\ne 5,9$, such that for each $i$ there
exists a prime dividing $p_i$ and not dividing $p_j$ for $j\ne i$ are linearly independent in the concordance group. 
\end{proposition}
\begin{proof}
Recall from Section 4.1 of \cite{MR1957829} (see also \cite{owens-strle:qhs4b}) that the $d$-invariants of the double-branched
cover of $K_{p,2}$ are
$$d(k)=\frac14-\frac{k^2}{2p}+\left\{\begin{array}{rl}
\frac14 & \text{\ if\ }\frac{p+1}2 +k \text{\ is even}\\
-\frac14 & \text{\ if\ }\frac{p+1}2 +k \text{\ is odd}
\end{array}\right.
$$
for $|k|\le \frac{p-1}2$, where $k=0$ corresponds to the spin structure. Then $d(0)\ne 0$ iff $p\equiv 3 \pmod 4$.
Assume now that $p\equiv 1 \pmod 4$, let $q$ be a prime dividing $p$ and write $p=qs$. Then we have (up to sign)
$$\mathcal{D}_q=2\sum_{k=1}^{(q-1)/2} d(ks)=\frac{q-1}4\left(1-\frac{s(q+1)}6\right)+
\left\{\begin{array}{ll}
\frac12 & \text{\ if\ } \frac{q-1}2 \text{\ is odd}\\
0 & \text{\ if\ }\frac{q-1}2 \text{\ is even}
\end{array}\right.\, .
$$
Now $\mathcal{D}_3=0$ iff $s=3$ which corresponds to Stevedore's knot, and $\mathcal{D}_5=0$ iff $s=1$ which
corresponds to the figure-$8$ knot. Note that $\mathcal{D}_3<0$ if $s>3$ and $\mathcal{D}_5<0$ if $s>1$. 
Finally $\mathcal{D}_q<0$ for $q\ge 7$ except for $q=7$ and $s=1$, which corresponds to $p=7$ 
(which is not congruent to $1$ modulo $4$).

To show linear independence suppose a knot $K=\#_{i=1}^n m_iK_{p_i,2}$ is slice. Let $q_i$ be a prime dividing $p_i$ and not
dividing any other $p_j$. If $q_i>7$, then $\mathcal{D}(K_{p_i,2})\ne 0$ and $\mathcal{D}(K_{p_j,2})= 0$ for $j\ne i$.
Hence it follows that $m_i=0$. If $q_i=3$ and $p_i>9$, or $q_i=5$ and $p_i>5$, or $q_i=7$ and $p_i>7$ the same conclusion
holds. If $p_i=3$ or $p_i=7$ then a direct computation shows that the $\mathcal D$ invariant is nonzero, leading to the
same conclusion.
\end{proof}

\section{Computing $\tau$ invariants for preimages of $2$-bridge knots}\label{sec:comp2}

To compute the $\tau$ invariants associated to the preimage, 
$\widetilde{K}_{p,q}$, of the $2$-bridge knot $K_{p,q} \subset S^3$ 
inside its double-branched cover $-L(p,q)$, we 
used a computer program written in Mathematica.  This program  implemented the 
combinatorial description of the knot Floer homology of the preimage of a $2$-bridge knot inside its double-branched cover 
presented in \cite{GT0610238}.  We 
summarize the results of that paper and add some minor improvements 
to the combinatorial description of the $({\bf A},{\bf M})$ bigrading using results from 
the more recent work in \cite{GT0610559}.

Recall that we can associate to $\widetilde{K}_{p,q} \subset -L(p,q)$ 
a compatible $4$-pointed, genus $1$ Heegaard diagram which is a {\it 
twisted toroidal grid diagram} consisting of two parallel curves of 
slope $0$ and two of slope $\frac{p}{q}$, partitioning the torus into 
$2pq$ cells.

More specifically, we identify the universal cover of the torus with 
the plane: $$T^2 := \mathbb{R}^2/\mathbb{Z}^2.$$

\begin{figure}\begin{center}\epsfig{file=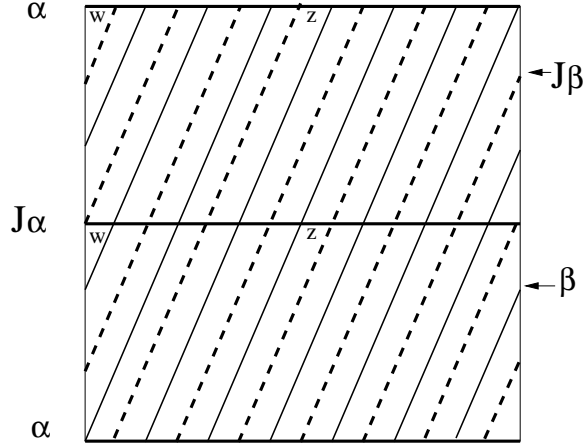, 
width=3in}\end{center}\caption{The twisted toroidal grid diagram 
which is a $4$-pointed genus $1$ Heegaard diagram for 
$\widetilde{K}_{7,3} \subset -L(7,3)$ (identify top-bottom, 
left-right, in the standard way).}\label{fig:K7_3Simple}\end{figure}

The two curves of slope $0$ on $T^2$ are the image in $T^2$ of the 
lines $y = 0$ and $y = \frac{1}{2}$ and the two curves of slope 
$\frac{p}{q}$ are the image in $T^2$ of the lines $y = \frac{p}{q}x$ 
and $y = \frac{p}{q}(x - \frac{1}{2})$.  We now identify the toroidal 
grid diagram with the fundamental domain $[0,1] \times [0,1] \subset 
\mathbb{R}^2$, \footnote{Note that by ``the image in $T^2$'' we mean 
the image of these lines in the quotient $\mathbb{R}^2/\mathbb{Z}^2$ 
and not the intersection of these lines with the chosen fundamental 
domain.}  and position our four basepoints at $$(\epsilon, 
1-\epsilon), (\frac{1}{2} + \epsilon, 1-\epsilon), (\epsilon, 
\frac{1}{2} - \epsilon), (\frac{1}{2} + \epsilon, \frac{1}{2} - 
\epsilon)$$ where $0<\epsilon<\min(\frac{1}{p}, \frac{1}{q})$.  See 
Figure \ref{fig:K7_3Simple} for the example of $K_{7,3}$.

We calculate $\tau_{\mathfrak{s}}$ for all $\mathfrak{s} \in 
\spinc(-L(p,q))$ by considering the chain complex

\begin{itemize}\item whose generators are indexed by bijections 
between the set of slope $0$ curves and the set of slope 
$\frac{p}{q}$ curves,\item whose differentials are given by counting 
parallelograms missing all $w$ basepoints, and\item with a filtration 
induced by the Alexander grading.\end{itemize}

By Lemma 2.1 in \cite{GT0610559}, this chain complex has the filtered 
chain homotopy type of $\cfhat(-L(p,q)) \otimes V$ where $V$ is 
the chain complex with two generators, one in bigrading ({\bf A}, 
{\bf M}) = $(0,0)$ and one in bigrading $(-1,-1)$, and no 
differentials.  In other words,

\begin{itemize}\item The associated graded complex of this chain 
complex is $$\hfkhat(-L(p,q),\widetilde{K}_{p,q}) \otimes V,$$ 
and\item the $E^{\infty}$ term of the spectral sequence associated to 
this filtered complex is $$\hfhat(-L(p,q)) \otimes V,$$  i.e., 
for all $\mathfrak{s} \in \spinc(-L(p,q))$, there are two 
generators, one in bigrading $(\tau_{\mathfrak{s}}, 
d_{\mathfrak{s}})$ and the other in bigrading 
$(\tau_{\mathfrak{s}}-1, d_{\mathfrak{s}}-1)$.\end{itemize}

\subsection{Enumerating Generators and Differentials} 
\begin{figure}\begin{center}\epsfig{file=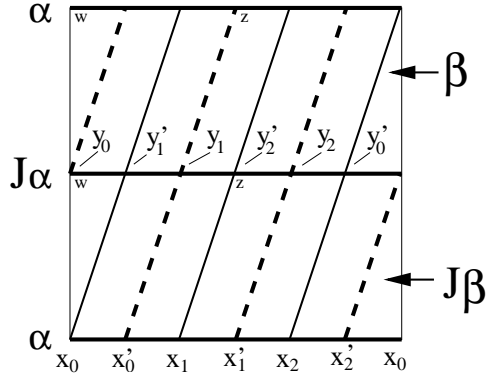, 
width=2.5in}\end{center}\caption{Labeling intersection points on the 
twisted toroidal grid diagram for $\widetilde{K}_{3,1} \subset 
-L(3,1)$.}\label{fig:K3_1Spinc}\end{figure}

We label the 
intersection points as described in \cite{GT0610238}.  Namely, 
intersection points with $\beta$ and $J(\beta)$ along the curve 
$\alpha$ are cyclically labeled $x_0, x_0', x_1, x_1', \ldots, 
x_{p-1}, x_{p-1}'$,  and intersection points with $\beta$ and 
$J(\beta)$ along the curve $J(\alpha)$ are labeled $y_0', y_0, y_1', 
y_1, \ldots, y_{p-1}', y_{p-1}$.  See Figure \ref{fig:K3_1Spinc}.

The generators of the chain complex are therefore pairs of the form 
$(x_i,y_j)$ and $(x_i',y_j')$ for $i,j \in \mathbb{Z}_p$.  This 
labeling convention is chosen so that the sum of indices constituting 
a generator specifies the \spinc structure in which it lives.  More 
precisely, $$(x_i,y_j), (x_i',y_j') \in \mathfrak{s}_{(i+j \mod 
p)},$$ where the subscript on $\mathfrak{s}$ denotes a particular 
affine identification of $\spinc(-L(p,q))$ with $H^2(-L(p,q),\zz) \cong 
\mathbb{Z}_p$, subject to the condition that the unique spin 
structure is identified with $0 \in \mathbb{Z}_p$.

There are $2p^2$ generators, $2p$ in each \spinc structure for $-L(p,q)$. 
See Table \ref{table:SpincGenerators} to see how the generators in 
Figure \ref{fig:K3_1Spinc} break up into \spinc structures.

\begin{table}\label{table:SpincGenerators}\begin{tabular}{|c|c|c|}\hline$\mathfrak{s}_0$ 
& $\mathfrak{s}_1$ & $\mathfrak{s}_2$\\\hline$(x_0,y_0)$   & 
$(x_0,y_1)$   & $(x_0,y_2)$\\$(x_0',y_0')$ & $(x_0',y_1')$ & 
$(x_0',y_2')$\\$(x_1,y_2)$   & $(x_1,y_0)$   & 
$(x_1,y_1)$\\$(x_1',y_2')$ & $(x_1',y_0')$ & 
$(x_1',y_1')$\\$(x_2,y_1)$   & $(x_2,y_2)$   & 
$(x_2,y_0)$\\$(x_2',y_1')$ & $(x_2',y_2')$ & 
$(x_2',y_0')$\\\hline\end{tabular}\caption{Splitting of generators in 
Figure \ref{fig:K3_1Spinc} according to \spinc 
structures}\end{table}

After enumerating the generators, we turn our attention to 
enumerating the differentials.  In other words, we wish to determine 
how many differentials (parallelograms) missing all $w$ basepoints 
connect ${\bf g_1}$ to ${\bf g_2}$, assuming that ${\bf g_1}$ and 
${\bf g_2}$ are two generators of our chain complex

\begin{itemize}\item living in the same \spinc structure and\item 
having relative homological (Maslov) grading\footnote{The computation 
of absolute $\mathbb{Q}$ homological gradings is addressed in the next subsection.} difference $1$, 
i.e., $${\bf M}({\bf g_1}) - {\bf M}({\bf g_2}) = 1.$$\end{itemize}

\begin{figure}\begin{center}\epsfig{file=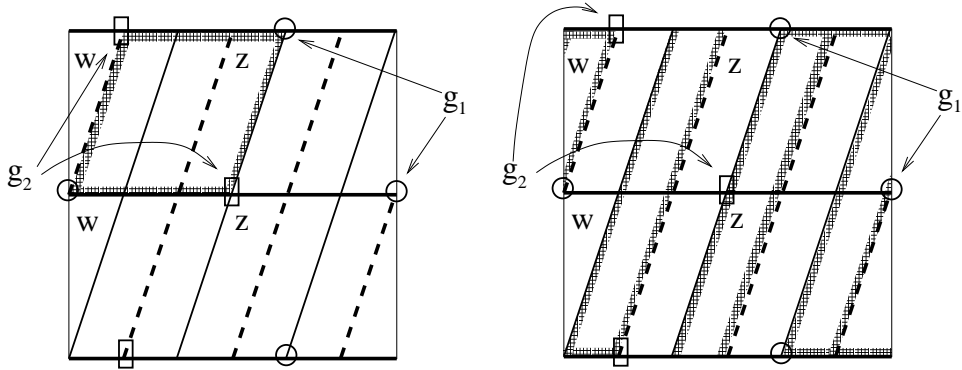, 
width=5in}\end{center}\caption{The two candidate differentials 
connecting ${\bf g_1}$ and ${\bf g_2}$ are shown.  Since the second 
parallelogram is non-imbedded and wraps around the torus several 
times, we have drawn its outline only.  Note that the first 
parallelogram misses all $w$ basepoints but the second does not. 
Therefore, the mod 2 multiplicity of ${\bf g_2}$ in the boundary of 
${\bf g_1}$ is 1.}\label{fig:Differentials1}\end{figure}

\begin{figure}\begin{center}\epsfig{file=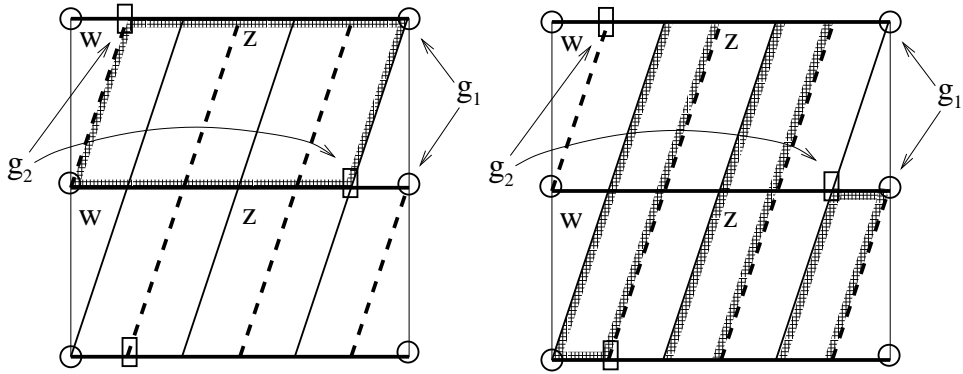, 
width=5in}\end{center}\caption{Here, both candidate parallelograms 
miss all $w$ basepoints.  Therefore, the mod 2 multiplicity of ${\bf 
g_2}$ in the boundary of ${\bf g_1}$ is 
0.}\label{fig:Differentials2}\end{figure}

Since there are always exactly $2$ parallelograms connecting any two 
generators in the same \spinc structure with relative Maslov 
grading $1$ (see Figures \ref{fig:Differentials1} and 
\ref{fig:Differentials2}), the only question is how many of these 
parallelograms (mod $2$) miss all $w$ basepoints.  Therefore, the 
multiplicity of ${\bf g_2}$ in the boundary of ${\bf g_1}$ is the 
reduction mod 2 of the number ($0,1,$ or $2$) of such parallelograms 
containing no $w$ basepoints.

\subsection{Computing Gradings}Denoting the set of generators of the 
chain complex by $\mathcal{G}$, the Alexander (filtration) grading is 
an assignment $${\bf A}:\mathcal{G} \rightarrow \mathbb{Z},$$ and the 
Maslov (homological) grading is an assignment $${\bf M}:\mathcal{G} 
\rightarrow \mathbb{Q}.$$

We begin by calculating the relative Maslov $\mathbb{Q}$-gradings of 
all generators.  We do this by lifting the pair $$(-L(p,q), 
\widetilde{K}_{p,q})$$ to its universal cover, 
$$(S^3,\widetilde{\widetilde{K}}_{p,q}),$$ calculating the relative 
Maslov gradings there using the easy formula proved in 
\cite{GT0610559}, then use Lee and Lipshitz's result in 
\cite{GT0608001} relating relative Maslov gradings under covers.  To 
nail down the absolute $\mathbb{Q}$-grading, we use the inductive 
formula for the correction terms in \cite{MR1957829} to pin down the 
absolute grading for one generator, thus pinning down the absolute 
grading for all generators.

The first step in this process is understanding how to construct a 
Heegaard diagram for the pair 
$$(S^3,\widetilde{\widetilde{K}}_{p,q})$$ from the Heegaard diagram 
for the pair $$(-L(p,q), \widetilde{K}_{p,q}).$$  The following lemma 
describes how to do this for any $2n$-pointed, grid number $n$, 
twisted toroidal grid diagram for a knot in a lens space (not just 
the grid number $2$ knots of interest here).

\begin{lemma}\label{lemma:CoverPicture}Let $T$ be a twisted toroidal 
grid diagram for $K$ in $L(p,q)$.  Form the universal cover, 
$\mathbb{R}^2$ of $T$, identifying $T$ with $$[0,1] \times [0,1] 
\subset \mathbb{R}^2,$$ the fundamental domain of the covering space 
action.  Let $Z$ be the lattice generated by the vectors $(1,0)$ and 
$(0,p)$.  Then $$\widetilde{T} = \mathbb{R}^2/Z$$ is a Heegaard 
diagram compatible with $\widetilde{K} \subset S^3$, where 
$\widetilde{K}$ is the preimage of $K$ under the covering space 
projection $\pi: S^3 \rightarrow L(p,q)$.\end{lemma}

{\flushleft \bf Proof of Lemma \ref{lemma:CoverPicture}:}  The 
original twisted toroidal grid diagram $T$ compatible with $K \subset 
L(p,q)$ corresponds to a handlebody decomposition of $L(p,q)$ with 
one solid handlebody, $Y_{0,1}$ formed by the union of $n$ $0$- and 
$1$-handles and the other solid handlebody $Y_{2,3}$ similarly formed 
by the union of $n$ $2$- and $3$-handles.  $\partial(Y_{0,1}) 
=\partial(-Y_{2,3}) = T.$

We will construct a Heegaard diagram for $\widetilde{K} \subset S^3$ 
by constructing a handlebody decomposition of the universal cover, 
$S^3$, compatible with this handlebody decomposition of $L(p,q)$ and 
the covering space action.  Namely,for every $h$ in the handlebody 
decomposition of $L(p,q)$, $\pi^{-1}(h) = 
\{\widetilde{h},a\widetilde{h},a^2\widetilde{h}, \ldots, 
a^{p-1}\widetilde{h}\}$ is a collection of handles in the handlebody 
decomposition for $S^3$, where $a$ is a generator of $\pi_1(L(p,q)) 
\cong \mathbb{Z}_p$, and $$\pi: S^3 \rightarrow L(p,q)$$ is the 
covering space map.  The attaching maps for the lifts of the handles 
are uniquely specified by the condition that they commute with the 
covering space projection.

Applying this procedure to the handlebody decomposition associated to 
$T$ corresponds to cutting $Y_{0,1}$ open along some meridian and 
gluing $p$ copies of the resulting $D^2 \times I$ together.  From the 
point of view of the boundary, $T$, this corresponds to stacking $p$ 
copies of $T$ on top of each other (when $T$ is identified with the 
fundamental domain $[0,1] \times [0,1]$ in $\mathbb{R}^2$).  This is 
precisely a description of $\widetilde{T} = \mathbb{R}^2/Z$.$\qed$

\begin{figure}\begin{center}\epsfig{file=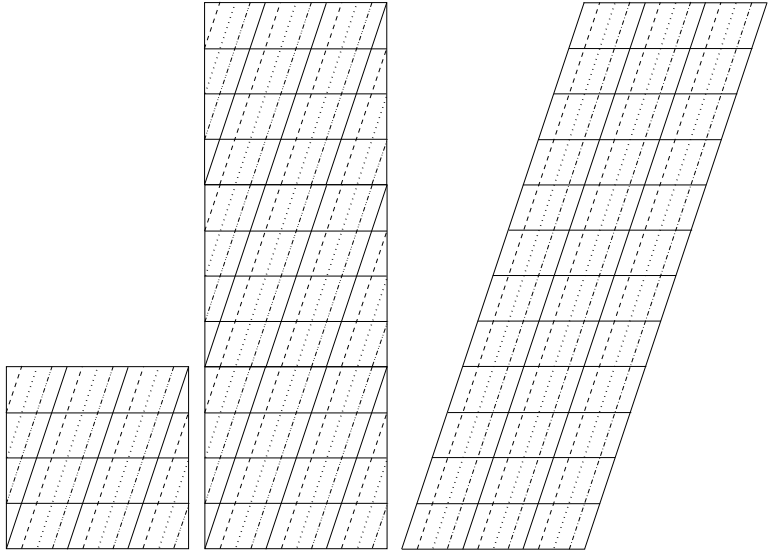, 
width=4in}\end{center}\caption{Constructing $\widetilde{T}$ from $T$ 
and adjusting the fundamental domain to identify $\widetilde{T}$ with 
a standard toroidal grid diagram for $\widetilde{K} \subset 
S^3$}\label{fig:SlantedDomain}\end{figure}

Note that $\widetilde{T}$ is just a standard (untwisted) toroidal 
grid diagram for the link $\widetilde{K}$ in $S^3$ (in the sense of 
\cite{GT0607691}).  With this in mind, it will be convenient for us 
to choose a slanted fundamental domain for $\widetilde{T} \subset 
\mathbb{R}^2,$ whose top and bottom edges are the same $\alpha$ curve 
and whose left and right edges are the same $\beta$ curve.  See 
Figure \ref{fig:SlantedDomain}.

We now recall the following fact, which is essentially Theorem 4.1 in 
\cite{GT0608001}.  See \cite{MR2222356} for a definition of 
$\widetilde{gr}$, the absolute $\mathbb{Q}$ homological grading.

\begin{theorem}\cite{GT0608001}Let ${\bf g_1}$ and ${\bf g_2}$ be 
generators in torsion \spinc structures in a Heegaard Floer chain 
complex, $\cfhat(Y)$, associated to a particular Heegaard 
decomposition $hd(Y)$ for $Y$.

Let $\pi:\widetilde{Y} \rightarrow Y$ be a degree $n$ covering map 
and $\widetilde{hd}(\widetilde{Y})$ the associated Heegaard 
decomposition of $\widetilde{Y}$ compatible with $\pi$.

Let $\widetilde{\bf g}_1$ and $\widetilde{\bf g}_2$ be the unique 
generators in $\cfhat(\widetilde{Y})$ with the property that 
$\pi(\widetilde{\bf g}_i) = {\bf g}_i$ for $i=1,2$.  Then 
$$\widetilde{gr}({\bf g}_1) - \widetilde{gr}({\bf g}_2) = 
\frac{1}{n}[\widetilde{gr}(\widetilde{{\bf g}}_1) - 
\widetilde{gr}(\widetilde{{\bf g}}_2)].$$\end{theorem}

This theorem allows us to compute the relative $\mathbb{Q}$ gradings 
between generators in the twisted toroidal grid diagram for 
$(-L(p,q), \widetilde{K}_{p,q})$ by lifting the generators to the 
$\pi$-compatible toroidal grid diagram for $(S^3, 
\widetilde{\widetilde{K}_{p,q}})$ and computing their relative 
gradings there.

Furthermore, in \cite{GT0610559}, an easy formula is given for 
determining the absolute $\mathbb{Q}$ grading for generators on a 
toroidal grid diagram for a link $K$ in $S^3$.  Namely, they define a 
function $$\mathcal{I}: {\bf S} \times {\bf S} \rightarrow 
\mathbb{Z}_+,$$ where ${\bf S}$ is the set of finite sets of  points 
on $\mathbb{R}^2$: If $A, B \in {\bf S}$, then $\mathcal{I}(A,B)$ is 
the number of pairs $(a_1,a_2), (b_1,b_2)$ for which $(a_1,a_2) \in 
A$, $(b_1,b_2) \in B$, $a_1 < a_2$ and $b_1 < b_2$.

Upon identification of the toroidal grid diagram with a fundamental 
domain in $\mathbb{R}^2$ with the property that the left edge is one 
of the $\beta$ curves and the bottom edge is one of the $\alpha$ 
curves, they then define a function:

$${\bf M}({\bf x}) = \mathcal{I}({\bf x},{\bf x}) + 
\mathcal{I}(\mathbb{O},\mathbb{O}) - \mathcal{I}({\bf x}, \mathbb{O}) 
- \mathcal{I}(\mathbb{O},{\bf x})+ 1$$where ${\bf x}$ is a generator 
of the chain complex, and ${\mathbb{O}}$ is the set of $w$ 
basepoints.  They then go on to prove that ${\bf M}$ (independent of 
the choice of identification of $\mathcal{T}$ with a fundamental 
domain on $\mathbb{R}^2$) is precisely $\widetilde{gr}$, the absolute 
homological grading, on $S^3$.

\cite{GT0610559} also describes how to obtain the Alexander grading 
for a generator by comparing its Maslov grading with respect to the 
$w$ basepoints with the Maslov grading with respect to the $z$ 
basepoints:
$${\bf A}({\bf x}) = \frac{1}{2}({\bf M}_w({\bf x}) - {\bf M}_z({\bf x})) - \frac{n-1}{2},$$
where $n$ is the grid number of $K$.  Although the formula is stated 
in \cite{GT0610559} only for a knot in $S^3$, it holds equally well 
for generators in a chain complex arising from a general balanced 
$2n$-pointed Heegaard decomposition of a $\qhs$. 
In fact, it is a direct consequence of one of the symmetries (see 
\cite{MR2065507}, \cite{GT0512286}) enjoyed by such a chain complex.

Namely, suppose $Y$ is a $\qhs$, and $K \subset 
Y$ is a nullhomogous knot.  Fix a balanced, $2n$-pointed Heegaard 
diagram for $Y$ compatible with $K$.  Switching the roles of the $w$ and $z$ 
basepoints on the same Heegaard diagram corresponds to reversing the 
orientation on $K$.  Furthermore, doing so induces a linear map on the chain complex sending a generator in bigrading $(i,d)$ to one in bigrading $(-i-(n-1),d-2i-(n-1))$.

Hence, we need only determine the absolute Maslov gradings with respect to the $w$ basepoints and then again with respect to the $z$ basepoints in order to compute all of the Alexander gradings.
\bibliography{taucover}\end{document}